\documentclass[reqno]{amsart}

\usepackage{graphicx}
\usepackage{amscd}
\usepackage{amsmath, amssymb}

\newtheorem{theorem}{Theorem}

\newtheorem{definition}[theorem]{Definition}
\newtheorem{example}[theorem]{Example}

\newtheorem{notation}[theorem]{Notation}

\newtheorem{proposition}[theorem]{Proposition}

\begin{document}

\title[ON THE FAITHFULNESS OF THE EXTENSION OF LAWRENCE-KRAMMER REPRESENTATION]{ON THE FAITHFULNESS OF THE EXTENSION OF LAWRENCE-KRAMMER REPRESENTATION OF THE GROUP OF CONJUGATING AUTOMORPHISMS $C_3$}

\author{Mohamad N. Nasser \and Mohammad N. Abdulrahim}

\address{Mohamad N. Nasser\\
         Department of Mathematics and Computer Science\\
         Beirut Arab University\\
         P.O. Box 11-5020, Beirut, Lebanon}
\email{m.nasser@bau.edu.lb}

\address{Mohammad N. Abdulrahim\\
         Department of Mathematics and Computer Science\\
         Beirut Arab University\\
         P.O. Box 11-5020, Beirut, Lebanon}
\email{mna@bau.edu.lb}

\maketitle

\begin{abstract}

Let $C_n$ be the group of conjugating automorphisms. We study the representation $\rho$ of $C_n$, an extension of Lawrence-Krammer representation of the braid group $B_n$, defined by Valerij G. Bardakov. As Bardakov proved that the representation $\rho$ is unfaithful for $n \geq 5$, the cases $n=3,4$ remain open. In our work, we make attempts towards the faithfulness of $\rho$ in the case $n=3$.

\end{abstract}

\medskip

\renewcommand{\thefootnote}{}

\vskip 0.1in

\section{Introduction} 

The braid group on $n$ strings, $B_n$, is the abstract group with generators $\sigma_1,\ldots,\sigma_{n-1}$ and a presentation as follows:
\begin{align*}
&\sigma_i\sigma_{i+1}\sigma_i = \sigma_{i+1}\sigma_i\sigma_{i+1} ,\hspace{0.5cm} i=1,2,\ldots,n-2,\\
&\sigma_i\sigma_j = \sigma_j\sigma_i , \hspace{2.1cm} |i-j|>2.
\end{align*}

Let $\mathbb{F}_n$ be a free group of $n$ generators $x_1,x_2,\ldots,x_n$. One of the generalizations of the braid group $B_n$ is the group of conjugating automorphisms $C_n$ [1]. Here $C_n$ is the subgroup of $Aut(\mathbb{F}_n)$ that satisfies for any $\phi \in C_n$, $\phi(x_i)={f_i}^{-1}x_{\Pi(i)}f_i$, where $\Pi$ is a permutation on $\{1,2,\ldots,n\}$ and $f_i=f_i(x_1,x_2,\ldots,x_n)$.

One of the most famous linear representations of $B_n$ is Lawrence-Krammer representation [4]. Braid groups are linear due to Lawrence-Krammer representations. It was shown that Lawrence-Krammer representations are faithful for all $n$ [2]. In [1], Bardakov uses Magnus representation defined in [3] to construct a linear representation $\rho: C_n \mapsto GL(V_n)$, where $V_n$ is a free module of dimension $n(n-1)/2$ with a basis $\{v_{i,j}\}, 1\leq i<j\leq n$. This representation is an extension of Lawrence-Krammer representation of $B_n$. It was shown that the representation $\rho$ is unfaithful for $n\geq 5$ [1]. However, the question of faithfulness of $\rho$ is still open for $n=3,4$.

We study, in section 3, the faithfulness of the representation $\rho$ for $n=3$. We prove that $\rho$ is unfaithful under some choices of $q$ (see Proposition 5). On the other hand, we prove that if $q^{6k}\neq 1$ for all $k \in \mathbb{Z}$ then the possible words in $\ker \rho$ are $A_1T^{s_1}A_2T^{s_2}\ldots A_{r-1}T^{s_{r-1}}A_rT^{s_r}$ and $T^{s_1}A_1T^{s_2}A_2\ldots T^{s_{r-1}}A_{r-1}T^{s_r}A_r$, where $T=\sigma_2\alpha_2\alpha_1, r\in \mathbb{N}, s_i \in \mathbb{Z}$ for all $1\leq i \leq r, \displaystyle \sum_{i=1}^{r}s_i =0, \displaystyle \sum_{i=1}^{r} length (A_i)$ is even and $A_i \in \{\alpha_1, \alpha_2, \alpha_1\alpha_2, \alpha_2\alpha_1, \alpha_1\alpha_2\alpha_1\}$ for all $1 \leq i \leq r$ (see Theorem 7). Moreover, we prove, under some conditions on $q$, that the words $A_1TA_2T\ldots A_{r-1}TA_rT^{1-r}$ and $T^{1-r}A_1TA_2\ldots TA_{r-1}TA_r$, where $r \in \mathbb{N}, \displaystyle \sum_{i=1}^{r} length (A_i)$ is even and $A_i \in \{\alpha_1, \alpha_2, \alpha_1\alpha_2, \alpha_2\alpha_1, \alpha_1\alpha_2\alpha_1\}$ for all $1 \leq i \leq r$ are not in $\ker \rho$ in the case $A_i=A_j$ for all $1\leq i,j \leq r$ (see Theorem 8). Also, we determine some conditions on $A_i$'s under which the words $A_1TA_2T\ldots A_{r-1}TA_rT^{1-r}$ and $T^{1-r}A_1TA_2T\ldots A_{r-1}TA_r$ are not in $\ker \rho$ (see Proposition 9, Theorem 11).

\section{Preliminaries} 

The group of conjugating automorphisms, $C_n$, is the subgroup of $Aut(\mathbb{F}_n)$ that satisfies for any $\phi \in C_n$, $\phi(x_i)={f_i}^{-1}x_{\Pi(i)}f_i$, where $\Pi$ is a permutation on $\{1,2,\ldots ,n\}$ and $f_i=f_i(x_1,x_2,\ldots ,x_n)$. Here $\mathbb{F}_n$ is the free group of $n$ generators $x_1,x_2,\ldots ,x_n$.

A.G. Savushkina [5] proved that the group of conjugating automorphisms $C_n$ is generated by automorphisms $\sigma_1,\sigma_2,\ldots ,\sigma_{n-1},\alpha_1,\alpha_2,\ldots ,\alpha_{n-1}$ of the free group $F_n$, where $\sigma_1,\sigma_2,\ldots ,\sigma_{n-1}$ generate the braid group $B_n$, and $\alpha_1,\alpha_2,\ldots ,\alpha_{n-1}$ generate the symmetric group $S_n$.

In [1], we see that the group $C_n$ is defined by the relations:
\begin{align*}
\sigma_i\sigma_{i+1}\sigma_i = \sigma_{i+1}\sigma_i\sigma_{i+1} ,\hspace{.3cm} &\text{for} \hspace{.3cm} i=1,2,\ldots ,n-2,\vspace{0.1cm}\\ 
\sigma_i\sigma_j = \sigma_j\sigma_i , \hspace{.3cm} &\text{for} \hspace{.3cm} |i-j|>2, \vspace{0.1cm}\\
\alpha^{2}_i=1, \hspace{.3cm} &\text{for} \hspace{.3cm} i=1,2,\ldots ,n-1,  \vspace{0.1cm}\\
\alpha_i\alpha_j=\alpha_j\alpha_i, \hspace{.3cm} &\text{for} \hspace{.3cm} |i-j|\geq 2,  \vspace{0.1cm}\\
\alpha_i\sigma_j=\sigma_j\alpha_i, \hspace{.3cm} &\text{for} \hspace{.3cm} |i-j|\geq 2, \vspace{0.1cm}\\
\sigma_i\alpha_{i+1}\alpha_i=\alpha_{i+1}\alpha_i\sigma_{i+1}, \hspace{.3cm} &\text{for} \hspace{.3cm} i=1,2,\ldots ,n-2, \vspace{0.1cm}\\
\sigma_{i+1}\sigma_{i}\alpha_{i+1}=\alpha_{i}\sigma_{i+1}\sigma_{i}, \hspace{.3cm} &\text{for} \hspace{.3cm} i=1,2,....,n-2.
\end{align*}

\begin{definition}
\text{[1]}
Let $V_n$ be a free module of dimension $n(n-1)/2$ and a basis $\{v_{i,j}\}, 1\leq i< j\leq n$ over the ring $\mathbb{Z}[q^{\pm1}]$ of Laurent polynomials in one variable. We introduce the representation $\rho:C_n \mapsto GL(V_n)$ by the actions of $\sigma_i's$ and $\alpha_i's$ on the basis $\{v_{i,j}\}$ as follows:
\begin{align*}
\left\{\begin{array}{l}
\sigma_i(v_{k,i})=(1-q)v_{k,i}+qv_{k,{i+1}}+q(q-1)v_{i,i+1},\\
\sigma_i(v_{k,{i+1}})=v_{k,i}$, \hspace{0.5cm} $k<i,\\
\sigma_i(v_{i,{i+1}})=q^2v_{i,{i+1}},\\
\sigma_i(v_{i,l})=q(q-1)v_{i,{i+1}}+(1-q)v_{i,l}+qv_{i+1,l},$ \hspace{0.5cm}$ i+1<l,\\
\sigma_i(v_{i+1,l})=v_{i,l},\\
\sigma_i(v_{k,l})=v_{k,l}$, \hspace{0.5cm} $\{ k,l\} \cap \{i,i+1\}=\emptyset, \\
\alpha_i(v_{k,i})=v_{k,{i+1}},\\
\alpha_i(v_{k,{i+1}})=v_{k,i}$, \hspace{0.5cm} $k<i,\\
\alpha_i(v_{i,{i+1}})=v_{i,{i+1}},\\
\alpha_i(v_{i,l})=v_{i+1,l},$ \hspace{0.5cm}$ i+1<l,\\
\alpha_i(v_{i+1,l})=v_{i,l},\\
\alpha_i(v_{k,l})=v_{k,l}$, \hspace{0.5cm} $\{ k,l\} \cap \{i,i+1\}=\emptyset. \\
\end{array}\right.
\end{align*}
\end{definition}
\vspace{0.2cm}
\begin{notation}
Let $x$, a word in $C_n$, be written as a product of powers of generators of $C_n$ and their inverses.  We denote the length of $x$ to be the sum of the absolute value of the powers of the generators and their inverses. For example, if $x=\sigma_1^{5}\alpha_2\sigma_2^{-2}\alpha_1^{-1}\sigma_1^{-2}$ then length($x$)= $|5|+|1|+|-2|+|-1|+|-2|=11$. 
\end{notation}

\section{The faithfulness of the representation $\rho$ for $n=3$} 

We know that Lawrence-Krammer representations of $B_n$ are faithful for all $n$ [2]. The representation $\rho$ is an extension to $C_n$ of Lawrence-Krammer representations of $B_n$ [1].  Bardakov proved that the extension $\rho$ is unfaithful for $n\geq 5$ [1]. The question of faithfulness of the representation $\rho$ is still open for $n=3,4$.

Now, we specialize $q$ to be a non zero complex number.

\begin{definition}

Consider the complex space $\mathbb{C}^3$, a free module of dimension $3$ with the canonical basis $\{e_1,e_2,e_3\}$ over the ring $\mathbb{Z}[q^{\pm1}]$ of Laurent polynomials in one variable. The representation $\rho:C_3 \mapsto GL(\mathbb{C}^3)$ is defined by the actions of $\sigma_1$, $\sigma_2$, $\alpha_1$ and $\alpha_2$ on the basis $\{e_1,e_2,e_3\}$ as follows:\vspace{0.3cm} \\
$$\sigma_1 \mapsto \left\{\begin{array}{l}
e_1\mapsto q^2e_1\\
e_2\mapsto q(q-1)e_1+(1-q)e_2+qe_3\\
e_3\mapsto e_2
\end{array}\right., \hspace{0.2cm} \sigma_2 \mapsto \left\{\begin{array}{l}
e_1\mapsto (1-q)e_1+qe_2+q(q-1)e_3\\
e_2\mapsto e_1\\
e_3\mapsto q^2e_3
\end{array}\right.,$$
\\
$$\alpha_1 \mapsto \left\{\begin{array}{l}
e_1\mapsto e_1\\
e_2\mapsto e_3\\
e_3\mapsto e_2
\end{array}\right.
\hspace{0.2cm} \text{and} \hspace{0.3cm} \alpha_2 \mapsto \left\{\begin{array}{l}
e_1\mapsto e_2\\
e_2\mapsto e_1\\
e_3\mapsto e_3
\end{array}\right..$$

\vspace{0.5cm}
In other words, for $n=3$, the representation $\rho$ is given by 
$$\rho:C_3 \mapsto GL(\mathbb{C}^3)$$
$$\rho(\sigma_1)=
\left( \begin{array}{@{}c@{}}
\begin{matrix}
   		q^2 & 0 & 0 \\
    	q(q-1) & 1-q & q \\
        0 & 1 & 0 \\
\end{matrix}
\end{array} \right), \hspace{0.5cm}
\rho(\sigma_2)=
\left( \begin{array}{@{}c@{}}
\begin{matrix}
   		1-q & q & q(q-1) \\
    	1 & 0 & 0 \\
        0 & 0 & q^2 \\
\end{matrix}
\end{array} \right),$$
\\
$$\rho(\alpha_1)=
\left( \begin{array}{@{}c@{}}
\begin{matrix}
   		1 & 0 & 0 \\
    	0 & 0 & 1 \\
        0 & 1 & 0 \\
\end{matrix}
\end{array} \right) \hspace{0.2cm}\text{and} \hspace{0.5cm}
\rho(\alpha_2)=
\left( \begin{array}{@{}c@{}}
\begin{matrix}
   		0 & 1 & 0 \\
    	1 & 0 & 0 \\
        0 & 0 & 1 \\
\end{matrix}
\end{array} \right).$$
\end{definition}

\vspace{0.5cm}

Let $T=\sigma_2\alpha_2\alpha_1$. We see that $\sigma_2=T\alpha_1\alpha_2$ and $\sigma_1=\alpha_2\alpha_1T\alpha_2\alpha_1$ and so $T, \alpha_1,$ and $\alpha_2$ generate $C_3$. Notice that $\rho(T)= \left( \begin{array}{@{}c@{}}
\begin{matrix}
   		q & q(q-1) & 1-q \\
    	0 & 0 & 1 \\
        0 & q^2 & 0 \\
\end{matrix}
\end{array} \right)$, and so $\rho(T^2)=q^2I_3 \in Z(GL(\mathbb{C}^3))$, the center of $GL(\mathbb{C}^3)$. So we get $\rho(T^{2k})=q^{2k}I_3$ and $\rho(T^{2k+1})=q^{2k}\rho(T)$ for all $k\in Z$.
\vspace{0.2cm}

\begin{proposition}
The words in $C_3$, which are written only as products of $\alpha_1$ and $\alpha_2$, are $\alpha_1$, $\alpha_2$, $\alpha_1\alpha_2$, $\alpha_2\alpha_1$, and $\alpha_1\alpha_2\alpha_1$.
\end{proposition}
\begin{proof}
We have $\alpha_i^{2}=1$ for $i=1,2$, which means that $\alpha_i^{-1}=\alpha_i$ for $i=1,2$. By direct computations, we see that $\alpha_1\alpha_2\alpha_1\alpha_2\alpha_1\alpha_2=1$. So any word in $C_3$, which is written only as products of $\alpha_1$ and $\alpha_2$, must be of length less than or equal to five.
\begin{itemize}
\item \underline{Length 1}: $\alpha_1$ and $\alpha_2$.
\item \underline{Length 2}: $\alpha_1\alpha_2$ and $\alpha_2\alpha_1$.
\item \underline{Length 3}: $\alpha_1\alpha_2\alpha_1$ and  $\alpha_2\alpha_1\alpha_2$. Since $\alpha_1\alpha_2\alpha_1\alpha_2\alpha_1\alpha_2=1$, it follows that $\alpha_1\alpha_2\alpha_1=\alpha_2\alpha_1\alpha_2$.
\item \underline{Length 4}: $\alpha_1\alpha_2\alpha_1\alpha_2$ and $\alpha_2\alpha_1\alpha_2\alpha_1$. Since $\alpha_1\alpha_2\alpha_1\alpha_2\alpha_1\alpha_2=1$, it follows that $\alpha_1\alpha_2\alpha_1\alpha_2=\alpha_2\alpha_1$ (a word of length 2) and $\alpha_2\alpha_1\alpha_2\alpha_1=\alpha_1\alpha_2$ (a word of length 2).
\item \underline{Length 5}: $\alpha_1\alpha_2\alpha_1\alpha_2\alpha_1$ and $\alpha_2\alpha_1\alpha_2\alpha_1\alpha_2$. Since $\alpha_1\alpha_2\alpha_1\alpha_2\alpha_1\alpha_2=1$, it follows that $\alpha_1\alpha_2\alpha_1\alpha_2\alpha_1=\alpha_2$ (a word of length 1) and $\alpha_2\alpha_1\alpha_2\alpha_1\alpha_2=\alpha_1$ (a word of length 1).
\end{itemize}
Therefore, the words in $C_3$, which are written only as products of $\alpha_1$ and $\alpha_2$, are $\alpha_1$, $\alpha_2$, $\alpha_1\alpha_2$, $\alpha_2\alpha_1$, and $\alpha_1\alpha_2\alpha_1$.
\end{proof}

We now show that the representation $\rho$ is unfaithful for some values of $q$.

\vspace{0.2cm}
Before we do that, we define the following sets. For $k \in \mathbb{Z}$, set \vspace{0.1cm}
\begin{itemize}
\item $P_k=\{q\in \mathbb{C^*}-\{4\}$, $q$ is a solution of $(2x)^{-2k}(x^2-2x-\sqrt{x-4}x^{3/2})^{2k}=1$ and $(2x)^{-2k}(x^2-2x+\sqrt{x-4}x^{3/2})^{2k}=1\}$,\vspace{0.1cm}
\item $R_k=\{q\in \mathbb{C^*}-\{1/4\}$, $q$ is a solution of $(2x)^{-2k}(1-2x-\sqrt{1-4x})^{2k}=1$ and $(2x)^{-2k}(1-2x+\sqrt{1-4x})^{2k}=1\}$,\vspace{0.1cm}
\item $S_k=\{q\in \mathbb{C^*}$, $q$ is a solution of $x^{-2k}(1-x)^{2k}=1\}$.
\end{itemize}

Notice that the sets $P_k, R_k,$ and $S_n$ are not empty sets for any even $k$ $\in \mathbb{Z}$ and any $n \in \mathbb{Z}$. We have:
\begin{itemize}
\item $2 \in P_k$ since $(2\times2)^{-2k}(2^2-2\times2 -\sqrt{2-4}2^{3/2})^{2k}=(4)^{-2k}(-\sqrt{-16})^{2k}= (4)^{-2k}(4i)^{2k}=i^{2k}=1$ because $k$ is even, and $(2\times2)^{-2k}(2^2-2\times2 +\sqrt{2-4}2^{3/2})^{2k}=(4)^{-2k}(\sqrt{-16})^{2k}= (4)^{-2k}(4i)^{2k}=i^{2k}=1$ because $k$ is even. Hence $P_k \neq \emptyset$ for any even $k$ $\in \mathbb{Z}$.
\item $\frac{1}{2} \in R_k$ since ($2\times \frac{1}{2})^{-2k}(1-2\times \frac{1}{2} -\sqrt{1-4\times \frac{1}{2}})^{2k}=(-\sqrt{-1})^{2k}=i^{2k}=1$ because $k$ is even, and ($2\times \frac{1}{2})^{-2k}(1-2\times \frac{1}{2}+\sqrt{1-4\times \frac{1}{2}})^{2k}=(\sqrt{-1})^{2k}=i^{2k}=1$ because $k$ is even. Hence $R_k \neq \emptyset$ for any even $k$ $\in \mathbb{Z}$.
\item $ \frac{1}{2} \in S_n$ since $(\frac{1}{2})^{-2n}(1-\frac{1}{2})^{2n}=1$. Hence $S_n \neq \emptyset$ for any $n \in \mathbb{Z}$.
\end{itemize}
\vspace{0.2cm}
\begin{proposition}
Suppose that $q^{2k}\neq 1$ for all $k \in \mathbb{Z}$. If there exists an even integer $m$ such that $q \in P_m\cup R_m$ or an integer $n$ such that $q\in S_n$ then $\rho$ is unfaithful.
\end{proposition}

\begin{proof}
Fix an even integer $m$ such that $q \in P_m\cup R_m$, then $q \in P_m$ or $q \in R_m$.\\
In the case $q \in P_m$, consider the word $x=(\alpha_2T)^{2m}T^{-2m}$. Suppose to get a contradiction that $x$ is a trivial word, then $(\alpha_2T)^{2m}T^{-2m}(e_1)=e_1$, and so $(\alpha_2T)^{2m}(e_1)=e_1$ since $T(e_1)=e_1$. This implies that $(\alpha_2T)^{2m+1}(e_1)=\alpha_2T(e_1)$, and so $(\alpha_2T)^{2m}(e_2)=e_2$ since $\alpha_2T(e_1)=e_2$. This also implies that $(\alpha_2T)^{2m+1}(e_2)=\alpha_2T(e_2)$, and so $(\alpha_2T)^{2m}(q^2e_3)=q^2e_3$ since $\alpha_2T(e_2)=q^2e_3$. Hence $(\alpha_2T)^{2m}(e_3)=e_3$ which implies that $(\alpha_2T)^{2m}$ is a trivial word and so $x=(\alpha_2T)^{2m}T^{-2m}=T^{-2m}$. Thus $T^{-2m}$ is a trivial word, which is a contradiction since $\rho(T^{-2m})=q^{-2m}I_3$ and $q^{-2m}\neq 1$. Therefore $x$ is not a trivial word. Now, we choose a certain basis to diagonalize the matrix $\rho(x)$ and so we get $\rho(x)= \vspace*{0.2cm}
\\
P\left( \begin{array}{@{}c@{}}
\begin{matrix}
   		1 & 0 & 0 \\
    	0 & (2q)^{-2m}(q^2-2q-\sqrt{q-4}q^{3/2})^{2m} & 0 \\
        0 & 0 & (2q)^{-2m}(q^2-2q+\sqrt{q-4}q^{3/2})^{2m} \\
\end{matrix}
\end{array} \right)P^{-1},$ \vspace*{0.2cm}
\\
where $P$ is the matrix of eigenvectors of $\rho(x)$. Since $q\in P_m$, it follows that $\rho(x)=I_3$ and so $x\in \ker \rho$. Hence $\rho$ is unfaithful.\\
In the case $q \in R_m$, consider the word $x=(\alpha_1\alpha_2\alpha_1T)^{2m}T^{-2m}$. Suppose to get a contradiction that $x$ is a trivial word, then $(\alpha_1\alpha_2\alpha_1T)^{2m}T^{-2m}(e_1)=e_1$, and so $(\alpha_1\alpha_2\alpha_1T)^{2m}(e_1)=e_1$ since $T(e_1)=e_1$. This implies that $(\alpha_1\alpha_2\alpha_1T)^{2m+1}(e_1)=\alpha_1\alpha_2\alpha_1T(e_1)$, and so $(\alpha_1\alpha_2\alpha_1T)^{2m}(e_3)=e_3$ since $\alpha_1\alpha_2\alpha_1T(e_1)=e_3$. This also implies that $(\alpha_1\alpha_2\alpha_1T)^{2m+1}(e_3)=\alpha_1\alpha_2\alpha_1T(e_3)$, and so $(\alpha_1\alpha_2\alpha_1T)^{2m}(q(q-1)e_1+qe_2+(1-q)e_3)=q(q-1)e_1+qe_2+(1-q)e_3$, hence $(\alpha_2T)^{2m}(e_2)=e_2$. This means that $(\alpha_1\alpha_2\alpha_1T)^{2m}$ is a trivial word, and so $x=(\alpha_1\alpha_2\alpha_1T)^{2m}T^{-2m}=T^{-2m}$. Thus $T^{-2m}$ is a trivial word, which is a contradiction since $\rho(T^{-2m})=q^{-2m}I_3$ and $q^{-2m}\neq 1$. Therefore $x$ is not a trivial word. Now, we choose a certain basis to diagonalize the matrix $\rho(x)$ and so we get $\rho(x)= \vspace*{0.2cm}\\ P\left( \begin{array}{@{}c@{}}
\begin{matrix}
   		(2q)^{-2m}(1-2q-\sqrt{1-4q})^{2m} & 0 & 0 \\
    	0 & (2q)^{-2m}(1-2q+\sqrt{1-4q})^{2m} & 0 \\
        0 & 0 & 1 \\
\end{matrix}
\end{array} \right)P^{-1},$\vspace*{0.2cm}
\\ where $P$ is the matrix of eigenvectors of $\rho(x)$. Since $q\in R_m$, it follows that $\rho(x)=I_3$ and so $x\in \ker \rho$. Hence $\rho$ is unfaithful.\\
In the same way, fix an integer $n$ such that $q \in S_n$ and consider the word $x=(T\alpha_1\alpha_2\alpha_1)^{2n}T^{-2n}$. Suppose to get a contradiction that $x$ is a trivial word, then $(T\alpha_1\alpha_2\alpha_1)^{2n}T^{-2n}(e_1)=e_1$, and so $(T\alpha_1\alpha_2\alpha_1)^{2n}(e_1)=e_1$ since $T(e_1)=e_1$. This implies that $(T\alpha_1\alpha_2\alpha_1)^{2n+1}(e_1)=T\alpha_1\alpha_2\alpha_1(e_1)$, and so $(T\alpha_1\alpha_2\alpha_1)^{2n}((1-q)e_1+qe_2+q(q-1)e_3)=(1-q)e_1+qe_2+q(q-1)e_3$, which gives that $(T\alpha_1\alpha_2\alpha_1)^{2n}(qe_2+q(q-1)e_3)=qe_2+q(q-1)e_3$ since $(T\alpha_1\alpha_2\alpha_1)^{2n}(e_1)=e_1$. On the other hand, $T\alpha_1\alpha_2\alpha_1(e_3)=e_1$ implies that $(T\alpha_1\alpha_2\alpha_1)^{2n}(e_3)=(T\alpha_1\alpha_2\alpha_1)^{2n-1}(e_1)$. But $(T\alpha_1\alpha_2\alpha_1)^{2n}(e_1)=e_1$ implies that $(T\alpha_1\alpha_2\alpha_1)^{2n-1}(e_1)=(T\alpha_1\alpha_2\alpha_1)^{-1}(e_1)$ and so $(T\alpha_1\alpha_2\alpha_1)^{2n}(e_3)=(T\alpha_1\alpha_2\alpha_1)^{-1}(e_1)=e_3$. Hence $(T\alpha_1\alpha_2\alpha_1)^{2n}(e_2)=e_2$ since $(T\alpha_1\alpha_2\alpha_1)^{2n}(qe_2+q(q-1)e_3)=qe_2+q(q-1)e_3$ and $(T\alpha_1\alpha_2\alpha_1)^{2n}(e_3)=e_3$. This means that $(T\alpha_1\alpha_2\alpha_1)^{2n}$ is a trivial word, and so $x=(T\alpha_1\alpha_2\alpha_1)^{2n}T^{-2n}=T^{-2n}$. Thus $T^{-2n}$ is a trivial word, which is a contradiction since $\rho(T^{-2n})=q^{-2n}I_3$ and $q^{-2n}\neq 1$. Therefore $x$ is not a trivial word. Now, $\rho(x)=\\  \left( \begin{array}{@{}c@{}}
\begin{matrix}
   		q^{-2n}(1-q)^{2n} & q(1-q^{-2n}(1-q)^{2n}) & 0 \\
    	0 & 1 & 0 \\
        0 & 0 & 1 \\
\end{matrix}
\end{array} \right)=I_3$ since $q \in S_n$. Hence $x \in \ker \rho$ and so $\rho$ is unfaithful.
\end{proof}
 
\begin{example}
\end{example} 
\noindent Take $k=2$. We have $q=2\in P_2$. Consider $x=(\alpha_2T)^{4}T^{-4}$. Notice that $x$ is not a trivial word, indeed\\
$x(e_1)=(\alpha_2T)^{4}T^{-4}(e_1) \\
\hspace*{.8cm}=\alpha_2T\alpha_2T\alpha_2T\alpha_2TT^{-4}(e_1)\\
\hspace*{.8cm}=\alpha_2T\alpha_2T\alpha_2T\alpha_2(e_1)\\
\hspace*{.8cm}=\alpha_2T\alpha_2T\alpha_2T(e_2)\\
\hspace*{.8cm}=\alpha_2T\alpha_2T\alpha_2(4e_3)\\
\hspace*{.8cm}=\alpha_2T\alpha_2T(4e_3)\\
\hspace*{.8cm}=\alpha_2T\alpha_2(-4e_1+8e_2+8e_3)\\
\hspace*{.8cm}=\alpha_2T(-4e_2+8e_1+8e_3)\\
\hspace*{.8cm}=\alpha_2(16e_2)\\
\hspace*{.8cm}=16e_1 \neq e_1.\\ $
By direct computations, we get $\rho(x)=I_3$ and so $x\in \ker\rho$. Hence $\rho$ in unfaithful.\\

We now determine conditions under which elements can possibly belong to $\ker \rho$.

\begin{theorem}
Suppose the $q^{6k}\neq 1$ for all $k\in \mathbb{Z}$, then the possible words in $\ker \rho$ are\\
(a) $A_1T^{s_1}A_2T^{s_2}\ldots A_{r-1}T^{s_{r-1}}A_rT^{s_r}$, \vspace{0.1cm}\\
(b) $T^{s_1}A_1T^{s_2}A_2\ldots T^{s_{r-1}}A_{r-1}T^{s_r}A_r$, \vspace{0.1cm}\\
where $r \in \mathbb{N}, s_i \in \mathbb{Z}$ for all $1 \leq i \leq r$, $\displaystyle \sum_{i=1}^{r}s_i =0, \displaystyle \sum_{i=1}^{r} length (A_i)$ is even and $A_i \in \{\alpha_1, \alpha_2, \alpha_1\alpha_2, \alpha_2\alpha_1, \alpha_1\alpha_2\alpha_1\}$ for all $1 \leq i \leq r$.
\end{theorem}

\begin{proof}
Let $x$ be a word in $C_3$ that is generated by $T$, $\alpha_1$ and $\alpha_2$. Then $x$ is either one of the following:\\
(i) $\alpha_1, \alpha_2, \alpha_1\alpha_2, \alpha_2\alpha_1, \alpha_1\alpha_2\alpha_1$,\vspace{0.1cm} \\
(ii) $T^k$, \vspace{0.1cm} \\
(iii) $A_1T^{s_1}A_2T^{s_2}\ldots A_{r-1}T^{s_{r-1}}A_rT^{s_r}$,\vspace{0.1cm} \\
(iv) $T^{s_1}A_1T^{s_2}A_2\ldots T^{s_{r-1}}A_{r-1}T^{s_r}A_r$,\vspace{0.1cm} \\
where $A_i \in \{\alpha_1, \alpha_2, \alpha_1\alpha_2, \alpha_2\alpha_1, \alpha_1\alpha_2\alpha_1\}$, $k \in \mathbb{N}$, $r \in \mathbb{N}$, and $s_i \in \mathbb{Z}$ for all $1 \leq i \leq r$.

\vspace{0.1cm}
Now we discuss each case separately: \vspace{0.1cm} \\
(i) If $x\in \{\alpha_1, \alpha_2, \alpha_1\alpha_2, \alpha_2\alpha_1, \alpha_1\alpha_2\alpha_1\}$, then we easily see that $x \notin \ker \rho$.\vspace{0.1cm} \\
(ii) If $x= T^k$, then $\det(\rho(T^k))=(-q)^{3k}\neq 1$ and so $T^k \notin \ker \rho$. \vspace{0.1cm} \\
(iii) If $x= A_1T^{s_1}A_2T^{s_2} \ldots A_{r-1}T^{s_{r-1}}A_rT^{s_r}$, then $det(\rho(x))=(-1)^{a_1+a_2+t-i}q^{3(t-i)}$, where $a_1=$ number of times $\alpha_1$ occurs in $x$, $a_2=$ number of times $\alpha_2$ occurs in $x$, $t=$ number of times $T$ occurs in $x$, and $i=$ number of times $T^{-1}$ occurs in $x$. If $t\neq i$, then $\det(\rho(x))\neq 1$ and so $x \notin \ker \rho$. Otherwise, if $t=i$ and $a_1+a_2$ is odd, then $\det(\rho(x))=-1$ and so $x \notin \ker \rho$. This leaves us with the only possibility of having elements in $\ker \rho$ with $\displaystyle \sum_{i=1}^{r}s_i =0$ and $\displaystyle \sum_{i=1}^{r} length (A_i)$ is even.\\
(iv) If $x= T^{s_1}A_1T^{s_2}A_2\ldots T^{s_{r-1}}A_{r-1}T^{s_r}A_r$, then this is similar to the previous case.
\end{proof}

Now, we determine few words that do not belong to $\ker \rho$. Before we do that, we define the following sets.
\begin{itemize}
\item $E=\{A_1T^{s_1}A_2T^{s_2}\ldots A_{r-1}T^{s_{r-1}}A_rT^{s_r}$ and $T^{s_1}A_1T^{s_2}A_2\ldots T^{s_{r-1}}A_{r-1}T^{s_r}A_r$, where $r \in \mathbb{N}, s_i \in \mathbb{Z}$ for all $1\leq i \leq r, \displaystyle \sum_{i=1}^{r}s_i =0, \displaystyle \sum_{i=1}^{r} length (A_i)$ is even and $A_i \in \{\alpha_1, \alpha_2, \alpha_1\alpha_2, \alpha_2\alpha_1, \alpha_1\alpha_2\alpha_1\}$ for all $1 \leq i \leq r\}$, the set of possible element in $\ker\rho$ in the case $q^{6k}\neq 1$ for all $k \in \mathbb{Z}$, \\
\item $E_1=\{x\in E, x=A_1TA_2T\ldots A_{r-1}TA_rT^{1-r}$ or $T^{1-r}A_1TA_2\ldots TA_{r-1}TA_r\}$.
\end{itemize}
Clearly we can see that $E_1\subset E$.
\vspace{.2cm}
\begin{theorem}
Suppose $q^{6k}\neq 1$ for all $k \in \mathbb{Z}$, $q \notin P_m\cup R_m$ for any even integer $m$ and $q\notin S_n$ for any  integer $n$. If $x \in E_1$ with $A_i=A_j$ for all $1\leq i,j \leq r$ then $x \notin  \ker \rho$. 
\end{theorem}

\begin{proof}
$x\in E_1$ implies that $x=A_1TA_2T\ldots A_{r-1}TA_rT^{1-r}$ or $x=T^{1-r}A_1TA_2\ldots TA_{r-1}TA_r$ where $r \in \mathbb{N}, \displaystyle \sum_{i=1}^{r} length (A_i)$ is even and $A_i \in \{\alpha_1, \alpha_2, \alpha_1\alpha_2, \alpha_2\alpha_1, \alpha_1\alpha_2\alpha_1\}$, for all $1 \leq i \leq r$. We consider the following five cases of $A_i$. \vspace{0.2cm}\\
(a) $A_i=\alpha_1$ for all $1\leq i\leq r$: \vspace{0.2cm} \\
(i) $x=(\alpha_1T)^rT^{-r}$. We have $\displaystyle \sum_{i=1}^{r} length (A_i)$ is even, and so $r=2k$, $k\in \mathbb{N}$. So $\rho(x)=\rho((\alpha_1T)^{2k}T^{-2k})= \left( \begin{array}{@{}c@{}}
\begin{matrix}
   		1 & (q^2-1) & q^{-2k}-1 \\
    	0 & q^{2k} & 0 \\
        0 & 0 & q^{-2k} \\
\end{matrix}
\end{array} \right).$ By our assumption, we have $q^{2k}\neq 1$. This implies that $\rho(x) \neq I_3$ and so $x\notin \ker \rho$. \vspace{0.2cm}\\
(ii) $x=T^{-r}(T\alpha_1)^r$. We have $\displaystyle \sum_{i=1}^{r} length (A_i)$ is even, and so $r=2k$, $k\in \mathbb{N}$. So $\rho(x)=\rho(T^{-2k}(T\alpha_1)^{2k})= \left( \begin{array}{@{}c@{}}
\begin{matrix}
   		1 & q^{-2k}-1 & (q^2-1) \\
    	0 & q^{-2k} & 0 \\
        0 & 0 & q^{2k} \\
\end{matrix}
\end{array} \right).$ By our assumption, we have $q^{2k}\neq 1$. This implies that $\rho(x) \neq I_3$ and so $x\notin \ker \rho$. \vspace{0.2cm}\\
(b) $A_i=\alpha_2$ for all $1\leq i\leq r$: \vspace{0.2cm}\\
(i) $x=(\alpha_2T)^rT^{-r}$. We have $\displaystyle \sum_{i=1}^{r} length (A_i)$ is even, and so $r=2k$, $k\in \mathbb{N}$. In the case $k$ is odd, we have $\rho(x)=\rho((\alpha_2T)^{2k}T^{-2k})= \left( \begin{array}{@{}c@{}}
\begin{matrix}
   		0 & * & * \\
    	* & * & * \\
        * & * & * \\
\end{matrix}
\end{array} \right)\neq I_3$. So $x \notin \ker\rho$. If $k$ is even and $q=4$ then $\rho(x)=\rho((\alpha_2T)^{2k}T^{-2k})= \left( \begin{array}{@{}c@{}}
\begin{matrix}
   		* & * & * \\
    	* & * & * \\
        * & * & 1-4k^2 \\
\end{matrix}
\end{array} \right),$ with $1-4k^2 \neq 1$, which means that $\rho(x) \neq I_3$. So $x \notin \ker\rho$. If $k$ is even and $q\neq 4$, then we choose a certain basis to diagonalize the matrix $\rho(x)$ and so we get $\rho(x)=\rho((\alpha_2T)^{2k}T^{-2k})= P\left( \begin{array}{@{}c@{}}
\begin{matrix}
   		1 & 0 & 0 \\
    	0 & (2q)^{-2k}(q^2-2q-\sqrt{q-4}q^{3/2})^{2k} & 0 \\
        0 & 0 & (2q)^{-2k}(q^2-2q+\sqrt{q-4}q^{3/2})^{2k} \\
\end{matrix}
\end{array} \right)P^{-1},$ where $P$ is the matrix of eigenvectors of $\rho(x)$. If $x \in \ker \rho$ then $q \in P_k$ for some even $k \in \mathbb{Z}$, which is a contradiction. So $x \notin \ker \rho$. \vspace{0.2cm}\\
(ii) $x=T^{-r}(T\alpha_2)^r$. We have $\displaystyle \sum_{i=1}^{r} length (A_i)$ is even, and so $r=2k$, $k\in \mathbb{N}$.\\
In the case $k$ is odd, we have $\rho(x)=\rho(T^{-2k}(T\alpha_2)^{2k})= \left( \begin{array}{@{}c@{}}
\begin{matrix}
   		0 & * & * \\
    	* & * & * \\
        * & * & * \\
\end{matrix}
\end{array} \right)\neq I_3$. So $x \notin \ker\rho$. If $k$ is even and $q=4$, then $\rho(x)=\rho(T^{-2k}(T\alpha_2)^{2k})= \left( \begin{array}{@{}c@{}}
\begin{matrix}
   		* & * & * \\
    	* & * & * \\
        * & * & 1-4k^2 \\
\end{matrix}
\end{array} \right),$ with $1-4k^2 \neq 1$, which means that $\rho(x) \neq I_3$. So $x \notin \ker\rho$. If $k$ is even and $q\neq 4$, then we choose a certain basis to diagonalize the matrix $\rho(x)$ and so we get $\rho(x)=\rho(T^{-2k}(T\alpha_2)^{2k})= P\left( \begin{array}{@{}c@{}}
\begin{matrix}
   		1 & 0 & 0 \\
    	0 & (2q)^{-2k}(q^2-2q-\sqrt{q-4}q^{3/2})^{2k} & 0 \\
        0 & 0 & (2q)^{-2k}(q^2-2q+\sqrt{q-4}q^{3/2})^{2k} \\
\end{matrix}
\end{array} \right)P^{-1},$ where $P$ is the matrix of eigenvectors of $\rho(x)$. If $x \in \ker \rho$ then $q \in P_k$ for some even $k \in \mathbb{Z}$, which is a contradiction. So $x \notin \ker \rho$.  \vspace{0.2cm} \\
(c) $A_i=\alpha_1\alpha_2$ for all $1\leq i\leq r$: \vspace{0.2cm}\\
(i) $x=(\alpha_1\alpha_2T)^rT^{-r}$. In the case $r=2k$, we have $\rho(x)=\rho((\alpha_1\alpha_2T)^{2k}T^{-2k})= \left( \begin{array}{@{}c@{}}
\begin{matrix}
   		* & * & * \\
    	* & q^{2k} & * \\
        * & * & * \\
\end{matrix}
\end{array} \right).$ By our assumption, we have $q^{2k}\neq 1$. This implies that $\rho(x) \neq I_3$ and so $x\notin \ker \rho$. In the case $r=2k+1$, we have $\rho(x)=\rho((\alpha_1\alpha_2T)^{2k+1}T^{-2k-1})= \left( \begin{array}{@{}c@{}}
\begin{matrix}
   		* & * & * \\
    	* & 0 & * \\
        * & * & * \\
\end{matrix}
\end{array} \right)\neq I_3,$ and so $x\notin \ker \rho$. \vspace{0.2cm}\\
(ii) $x=T^{-r}(T\alpha_1\alpha_2)^r$. In the case $r=2k$, we have $\rho(x)=\rho(T^{-2k}(T\alpha_1\alpha_2)^{2k})= \left( \begin{array}{@{}c@{}}
\begin{matrix}
   		* & * & * \\
    	* & * & * \\
        * & * & q^{2k} \\
\end{matrix}
\end{array} \right).$ By our assumption, we have $q^{2k}\neq 1$. This implies that $\rho(x) \neq I_3$ and so $x\notin \ker \rho$. In the case $r=2k+1$, we have $\rho(x)=\rho(T^{-2k-1}(T\alpha_1\alpha_2)^{2k+1})= \left( \begin{array}{@{}c@{}}
\begin{matrix}
   		* & * & * \\
    	* & * & 1 \\
        * & * & * \\
\end{matrix}
\end{array} \right)\neq I_3,$ and so $x\notin \ker \rho$.\vspace{0.2cm} \\
(d) $A_i=\alpha_2\alpha_1$ for all $1\leq i\leq r$: \vspace{0.2cm}\\
(i) $x=(\alpha_2\alpha_1T)^rT^{-r}$. In the case $r=2k$, we have $\rho(x)=\rho((\alpha_2\alpha_1T)^{2k}T^{-2k})= \left( \begin{array}{@{}c@{}}
\begin{matrix}
   		* & * & * \\
    	* & * & * \\
        * & * & q^{-2k} \\
\end{matrix}
\end{array} \right).$ By our assumption, we have $q^{-2k}\neq 1$. This implies that $\rho(x) \neq I_3$ and so $x\notin \ker \rho$. In the case $r=2k+1$, we have $\rho(x)=\rho((\alpha_2\alpha_1T)^{2k+1}T^{-2k-1})= \left( \begin{array}{@{}c@{}}
\begin{matrix}
   		* & * & * \\
    	* & * & * \\
        * & * & 0 \\
\end{matrix}
\end{array} \right)\neq I_3,$ and so $x\notin \ker \rho$. \vspace{0.2cm}\\
(ii) $x=T^{-r}(T\alpha_2\alpha_1)^r$. In the case $r=2k$, we have $\rho(x)=\rho(T^{-2k}(T\alpha_2\alpha_1)^{2k})= \left( \begin{array}{@{}c@{}}
\begin{matrix}
   		* & * & * \\
    	* & q^{-2k} & * \\
        * & * & * \\
\end{matrix}
\end{array} \right).$ By our assumption, we have $q^{-2k}\neq 1$. This implies that $\rho(x) \neq I_3$ and so $x\notin \ker \rho$. In the case $r=2k+1$, we have $\rho(x)=\rho(T^{-2k-1}(T\alpha_2\alpha_1)^{2k+1})= \left( \begin{array}{@{}c@{}}
\begin{matrix}
   		* & * & * \\
    	* & * & * \\
        * & * & 0 \\
\end{matrix}
\end{array} \right)\neq I_3$ and so $x\notin \ker \rho$. \vspace{0.2cm}\\
(e) $A_i=\alpha_1\alpha_2\alpha_1$ for all $1\leq i\leq r$: \vspace{0.2cm}\\
(i) $x=(\alpha_1\alpha_2\alpha_1T)^rT^{-r}$. We have $\displaystyle \sum_{i=1}^{r} length (A_i)$ is even and so $r=2k$, $k\in \mathbb{N}$. In the case $k$ is odd, we have $\rho(x)=\rho((\alpha_1\alpha_2\alpha_1T)^{2k}T^{-2k})= \left( \begin{array}{@{}c@{}}
\begin{matrix}
   		0 & * & * \\
    	* & * & * \\
        * & * & * \\
\end{matrix}
\end{array} \right)\neq I_3$. So $x \notin \ker\rho$. If $k$ is even and $q=1/4$ then $\rho(x)=\rho((\alpha_1\alpha_2\alpha_1T)^{2k}T^{-2k})= \left( \begin{array}{@{}c@{}}
\begin{matrix}
   		* & * & * \\
    	* & * & * \\
        k(2k+1) & * & * \\
\end{matrix}
\end{array} \right),$ with $k(2k+1) \neq 0$, which means that $\rho(x) \neq I_3$. So $x \notin \ker\rho$. If $k$ is even and $q\neq 1/4$, then we choose a certain basis to diagonalize the matrix $\rho(x)$ and so we get $\rho(x)=\rho((\alpha_1\alpha_2\alpha_1T)^{2k}T^{-2k})= \vspace{.1cm}\\ P\left( \begin{array}{@{}c@{}}
\begin{matrix}
   		(2q)^{-2k}(1-2q-\sqrt{1-4q})^{2k} & 0 & 0 \\
    	0 & (2q)^{-2k}(1-2q+\sqrt{1-4q})^{2k} & 0 \\
        0 & 0 & 1 \\
\end{matrix}
\end{array} \right)P^{-1},$ where $P$ is the matrix of eigenvectors of $\rho(x)$. If $x \in \ker \rho$ then $q \in R_k$ for some even $k \in \mathbb{Z}$, which is a contradiction. So $x \notin \ker \rho$. \vspace{0.2cm}\\
(ii) $x=T^{-r}(T\alpha_1\alpha_2\alpha_1)^r$. We have $\displaystyle \sum_{i=1}^rlength(A_i)$ is even, and so $r=2k$, $k\in \mathbb{N}$. So $\rho(x)=\rho(T^{-2k}(T\alpha_1\alpha_2\alpha_1)^{2k})= \left( \begin{array}{@{}c@{}}
\begin{matrix}
   		q^{-2k}(1-q)^{2k} & q(1-q^{-2k}(1-q)^{2k}) & 0 \\
    	0 & 1 & 0 \\
        0 & 0 & 1 \\
\end{matrix}
\end{array} \right) \neq I_3$ since $q \notin S_m$ for all $m \in \mathbb{Z}$. So $x \notin \ker \rho$
\end{proof}
Next, we determine further conditions under which we eliminate few words from belonging to the kernel.
\begin{proposition}
Let $r\geq 3$ be an odd integer and let $x \in E_1$. Suppose one of the following holds true. \vspace{0.1cm}\\
(a) $A_rA_1=A_2=A_3=\ldots =A_{r-1}$. \vspace{0.1cm}\\
(b) $A_rA_1=A_{r-1}A_2=\ldots=A_{\frac{r+3}{2}}A_{\frac{r-1}{2}}=1$.\\
(c) There exists $i \in \{0,1,2,\ldots ,\frac{r-5}{2}\}$ such that $A_rA_1=A_{r-1}A_2=\ldots =A_{r-i}A_{i+1}=1$ and $A_{r-i-1}A_{i+2}=A_{i+3}=A_{i+4}=\ldots =A_{r-i-2}$.\vspace{0.2cm}\\
Then there exists $w \in C_3$ such that $\rho(w^{-1}xw)=\rho(\underbrace{AT}_1AT\ldots \underbrace{AT}_sT^{-s})$, where $A \in \{\alpha_1, \alpha_2, \alpha_1\alpha_2, \alpha_2\alpha_1, \alpha_1\alpha_2\alpha_1\}$ and $s\in \mathbb{N}$.
\end{proposition}
\begin{proof}
Without loss of generality, we assume that $x=A_1TA_2T\ldots A_{r-1}TA_rT^{1-r}.$\\
If (a) holds true, then we take $w=A_1T$. We have \\
$\rho(w^{-1}xw)=\rho((A_1T)^{-1}A_1TA_2T\ldots A_{r-1}TA_rT^{-r+1}(A_1T))$ \\
\hspace*{1.55cm} $=\rho(A_2T\ldots A_{r-1}TA_rT^{-r+1}A_1T)$ \\
\hspace*{1.55cm} $=\rho(A_2)\rho(T)\ldots \rho(A_{r-1})\rho(T)\rho(A_r)\rho(T^{-r+1})\rho(A_1)\rho(T)$ \\
\hspace*{1.55cm} $=\rho(A_2)\rho(T)\ldots \rho(A_{r-1})\rho(T)\rho(A_r)\rho(A_1)\rho(T^{-r+1})\rho(T)$ \\
\hspace*{1.55cm} $= \rho(A_2)\rho(T)\ldots \rho(A_{r-1})\rho(T)\rho(A_r)\rho(A_1)\rho(T^{-r+1}T)$,\\
\hspace*{1.55cm} $= \rho(A_2)\rho(T)\ldots \rho(A_{r-1})\rho(T)\rho(A_r)\rho(A_1)\rho(T^{-r+2})$,\\
\hspace*{1.55cm} $= \rho(A_2T\ldots A_{r-1}TA_rA_1T^{-r+2})$,\\
where $A_rA_1=A_2=\ldots=A_{r-1}$. \vspace{0.1cm}\\
If (b) holds true, then we take $w=A_1TA_2TA_3T\ldots A_{\frac{r-1}{2}}T$. We have \\
$\rho(w^{-1}xw)=\rho((A_1TA_2TA_3T\ldots A_{\frac{r-1}{2}}T)^{-1}A_1TA_2T\ldots A_{r-2}TA_{r-1}TA_rTT^{-r}(A_1TA_2TA_3T\ldots A_{\frac{r-1}{2}}T))$ \\
\hspace*{1.55cm} $=\rho((A_1T\ldots A_{\frac{r-1}{2}}T)^{-1}(A_1T\ldots A_{\frac{r-1}{2}}T)(A_{\frac{r+1}{2}}TA_{\frac{r+3}{2}}T\ldots A_{r-2}TA_{r-1}TA_rTT^{-r})(A_1T\ldots A_{\frac{r-1}{2}}T))$ \\
\hspace*{1.55cm} $=\rho(A_{\frac{r+1}{2}}TA_{\frac{r+3}{2}}T\ldots A_{r-2}TA_{r-1}TA_rT^{-r+1}A_1TA_2TA_3T\ldots A_{\frac{r-1}{2}}T)$ \\
\hspace*{1.55cm} $=\rho(A_{\frac{r+1}{2}})\rho(T)\rho(A_{\frac{r+3}{2}})\rho(T)\ldots \rho(A_{r-2})\rho(T)\rho(A_{r-1})\rho(T)\rho(A_r)\rho(T^{-r+1})\rho(A_1)\rho(T)\\
\hspace*{2cm} \rho(A_2)\rho(T)\rho(A_3)\rho(T)\ldots \rho(A_{\frac{r-1}{2}})\rho(T)$ \\
\hspace*{1.55cm} $=\rho(A_{\frac{r+1}{2}})\rho(T)\rho(A_{\frac{r+3}{2}})\rho(T)\ldots \rho(A_{r-2})\rho(T)\rho(A_{r-1})\rho(T)\rho(A_r)\rho(A_1)\rho(T^{-r+1})\rho(T)\\
\hspace*{2cm} \rho(A_2)\rho(T)\rho(A_3)\rho(T)\ldots \rho(A_{\frac{r-1}{2}})\rho(T)$ \\
\hspace*{1.55cm} $=\rho(A_{\frac{r+1}{2}})\rho(T)\rho(A_{\frac{r+3}{2}})\rho(T)\ldots \rho(A_{r-2})\rho(T)\rho(A_{r-1})\rho(T)\rho(A_rA_1)\rho(T^{-r+1})\rho(T)\\
\hspace*{2cm} \rho(A_2)\rho(T)\rho(A_3)\rho(T)\ldots \rho(A_{\frac{r-1}{2}})\rho(T)$ \\
\hspace*{1.55cm} $=\rho(A_{\frac{r+1}{2}})\rho(T)\rho(A_{\frac{r+3}{2}})\rho(T)\ldots \rho(A_{r-2})\rho(T)\rho(A_{r-1})\rho(T)\rho(T^{-r+1})\rho(T)\\
\hspace*{2cm}\rho(A_2)\rho(T)\rho(A_3)\rho(T)\ldots \rho(A_{\frac{r-1}{2}})\rho(T)$ \\
\hspace*{1.55cm} $=\rho(A_{\frac{r+1}{2}})\rho(T)\rho(A_{\frac{r+3}{2}})\rho(T)\ldots \rho(A_{r-2})\rho(T)\rho(A_{r-1})\rho(T^{-r+3})\rho(A_2)\rho(T)\\
\hspace*{2cm}\rho(A_3)\rho(T)\ldots \rho(A_{\frac{r-1}{2}})\rho(T)$ \\
\hspace*{1.55cm} $=\rho(A_{\frac{r+1}{2}}T)\rho(A_{\frac{r+3}{2}})\rho(T)\ldots \rho(A_{r-2})\rho(T)\rho(A_{r-1})\rho(A_2)\rho(T^{-r+3})\rho(T)\\
\hspace*{2cm}\rho(A_3)\rho(T)\ldots \rho(A_{\frac{r-1}{2}})\rho(T)$ \\
\hspace*{1.55cm} $=\rho(A_{\frac{r+1}{2}}T)\rho(A_{\frac{r+3}{2}})\rho(T)\ldots \rho(A_{r-2})\rho(T)\rho(A_{r-1}A_2)\rho(T^{-r+3})\rho(T)\\
\hspace*{2cm}\rho(A_3)\rho(T)\ldots \rho(A_{\frac{r-1}{2}})\rho(T)$ \\
\hspace*{1.55cm} $=\rho(A_{\frac{r+1}{2}})\rho(T)\rho(A_{\frac{r+3}{2}})\rho(T)\ldots \rho(A_{r-2})\rho(T)\rho(T^{-r+3})\rho(T)\rho(A_3)\rho(T)\ldots \rho(A_{\frac{r-1}{2}})\rho(T)$ \\
\hspace*{1.55cm} $=\rho(A_{\frac{r+1}{2}})\rho(T)\rho(A_{\frac{r+3}{2}})\rho(T)\ldots \rho(A_{r-2})\rho(T^{-r+5})\rho(A_3)\rho(T)\ldots \rho(A_{\frac{r-1}{2}})\rho(T)$ \\
\hspace*{1.65cm} \vdots \\
\hspace*{1.55cm} $=\rho(A_{\frac{r+1}{2}})\rho(T)\rho(A_{\frac{r+3}{2}})\rho(T^{-2})\rho(A_{\frac{r-1}{2}})\rho(T)$ \\
\hspace*{1.55cm} $=\rho(A_{\frac{r+1}{2}})\rho(T)\rho(A_{\frac{r+3}{2}})\rho(A_{\frac{r-1}{2}})\rho(T^{-2})\rho(T)$ \\
\hspace*{1.55cm} $=\rho(A_{\frac{r+1}{2}})\rho(T)\rho(A_{\frac{r+3}{2}}A_{\frac{r-1}{2}})\rho(T^{-2})\rho(T)$ \\
\hspace*{1.55cm} $=\rho(A_{\frac{r+1}{2}})\rho(T)\rho(T^{-2})\rho(T)$ \\
\hspace*{1.55cm} $=\rho(A_{\frac{r+1}{2}})$. \\
If (c) holds true, then there exists $i \in \{0,1,2,\ldots ,\frac{r-5}{2}\}$ such that $A_rA_1=A_{r-1}A_2=\ldots =A_{r-i}A_{i+1}=1$ and $A_{r-i-1}A_{i+2}=A_{i+3}=A_{i+4}=\ldots =A_{r-i-2}$. We take $w_i=A_1TA_2T\ldots A_{i+2}T$. We have\\
$\rho(w_i^{-1}xw_i)=\rho((A_1TA_2TA_3T\ldots A_{i+2}T)^{-1}A_1TA_2T\ldots A_{r-2}TA_{r-1}TA_rTT^{-r}(A_1TA_2TA_3T\ldots A_{i+2}T))$ \\
\hspace*{1.65cm} $=\rho((A_1T\ldots A_{i+2}T)^{-1}(A_1T\ldots A_{i+2}T)(A_{i+3}TA_{i+4}T\ldots A_{r-2}TA_{r-1}TA_rTT^{-r}(A_1T\ldots A_{i+2}T))$ \\
\hspace*{1.6cm} $=\rho(A_{i+3}TA_{i+4}T\ldots A_{r-2}TA_{r-1}TA_rT^{-r+1}A_1TA_2TA_3T\ldots A_{i+2}T)$ \\
\hspace*{1.65cm} $=\rho(A_{i+3})\rho(T)\rho(A_{i+4})\rho(T)\ldots \rho(A_{r-2})\rho(T)\rho(A_{r-1})\rho(T)\rho(A_r)\rho(T^{-r+1})\rho(A_1)\rho(T)\\
\hspace*{2.1cm}\rho(A_2)\rho(T)\rho(A_3)\rho(T)\ldots \rho(A_{i+2})\rho(T)$ \\
\hspace*{1.65cm} $=\rho(A_{i+3})\rho(T)\rho(A_{i+4})\rho(T)\ldots \rho(A_{r-2})\rho(T)\rho(A_{r-1})\rho(T)\rho(A_r)\rho(A_1)\rho(T^{-r+1})\rho(T)\\
\hspace*{2.1cm}\rho(A_2)\rho(T)\rho(A_3)\rho(T)\ldots \rho(A_{i+2})\rho(T)$ \\
\hspace*{1.65cm} $=\rho(A_{i+3})\rho(T)\rho(A_{i+4})\rho(T)\ldots \rho(A_{r-2})\rho(T)\rho(A_{r-1})\rho(T)\rho(A_rA_1)\rho(T^{-r+1})\rho(T)\\
\hspace*{2.1cm}\rho(A_2)\rho(T)\rho(A_3)\rho(T)\ldots \rho(A_{i+2})\rho(T)$ \\
\hspace*{1.65cm} $=\rho(A_{i+3})\rho(T)\rho(A_{i+4})\rho(T)\ldots \rho(A_{r-2})\rho(T)\rho(A_{r-1})\rho(T)\rho(T^{-r+1})\rho(T)\\
\hspace*{2.1cm}\rho(A_2)\rho(T)\rho(A_3)\rho(T)\ldots \rho(A_{i+2})\rho(T)$ \\
\hspace*{1.65cm} $=\rho(A_{i+3})\rho(T)\rho(A_{i+4})\rho(T)\ldots \rho(A_{r-2})\rho(T)\rho(A_{r-1})\rho(T^{-r+3})\rho(A_2)\rho(T)\\
\hspace*{2.1cm}\rho(A_3)\rho(T)\ldots \rho(A_{i+2})\rho(T)$ \\
\hspace*{1.65cm} $=\rho(A_{i+3})\rho(T)\rho(A_{i+4})\rho(T)\ldots \rho(A_{r-2})\rho(T)\rho(A_{r-1})\rho(A_2)\rho(T^{-r+3})\rho(T)$\\
\hspace*{1.65cm} $=\rho(A_{i+3})\rho(T)\rho(A_{i+4})\rho(T)\ldots \rho(A_{r-2})\rho(T)\rho(A_{r-1}A_2)\rho(T^{-r+3})\rho(T)\\
\hspace*{2.1cm}\rho(A_3)\rho(T)\ldots \rho(A_{i+2})\rho(T)$ \\
\hspace*{1.65cm} $=\rho(A_{i+3})\rho(T)\rho(A_{i+4})\rho(T)\ldots \rho(A_{r-2})\rho(T)\rho(T^{-r+3})\rho(T)\rho(A_3)\rho(T)\ldots \rho(A_{i+2})\rho(T)$ \\
\hspace*{1.65cm} $=\rho(A_{i+3})\rho(T)\rho(A_{i+4})\rho(T)\ldots \rho(A_{r-2})\rho(T^{-r+5})\rho(A_3)\rho(T)\ldots \rho(A_{i+2})\rho(T)$ \\
\hspace*{1.7cm} \vdots \\
\hspace*{1.65cm} $=\rho(A_{i+3})\rho(T)\rho(A_{i+4})\rho(T)\ldots \rho(A_{r-i-2})\rho(T)\rho(A_{r-i-1})\rho(T)\rho(A_{r-i})\rho(T^{-r+2i+1})\rho(A_{i+1})\\
\hspace*{2cm}\rho(T)\rho(A_{i+2})\rho(T)$\\
\hspace*{1.65cm} $=\rho(A_{i+3})\rho(T)\rho(A_{i+4})\rho(T)\ldots \rho(A_{r-i-2})\rho(T)\rho(A_{r-i-1})\rho(T)\rho(A_{r-i})\rho(A_{i+1})\rho(T^{-r+2i+1})\\
\hspace*{2.1cm}\rho(T)\rho(A_{i+2})\rho(T)$\\
\hspace*{1.65cm} $=\rho(A_{i+3})\rho(T)\rho(A_{i+4})\rho(T)\ldots \rho(A_{r-i-2})\rho(T)\rho(A_{r-i-1})\rho(T)\rho(A_{r-i}A_{i+1})\rho(T^{-r+2i+1})\\
\hspace*{2.1cm}\rho(T)\rho(A_{i+2})\rho(T)$\\
\hspace*{1.65cm} $=\rho(A_{i+3})\rho(T)\rho(A_{i+4})\rho(T)\ldots \rho(A_{r-i-2})\rho(T)\rho(A_{r-i-1})\rho(T)\rho(T^{-r+2i+1})\rho(T)\rho(A_{i+2})\rho(T)$\\
\hspace*{1.65cm} $=\rho(A_{i+3})\rho(T)\rho(A_{i+4})\rho(T)\ldots \rho(A_{r-i-2})\rho(T)\rho(A_{r-i-1})\rho(T^{-r+2i+3})\rho(A_{i+2})\rho(T)$\\
\hspace*{1.65cm} $=\rho(A_{i+3})\rho(T)\rho(A_{i+4})\rho(T)\ldots \rho(A_{r-i-2})\rho(T)\rho(A_{r-i-1})\rho(A_{i+2})\rho(T^{-r+2i+3})\rho(T)$\\
\hspace*{1.65cm} $=\rho(A_{i+3})\rho(T)\rho(A_{i+4})\rho(T)\ldots \rho(A_{r-i-2})\rho(T)\rho(A_{r-i-1}A_{i+2})\rho(T^{-r+2i+3})\rho(T)$\\
\hspace*{1.65cm} $= \rho(A_{i+3}TA_{i+4}T\ldots A_{r-i-2}TA_{r-i-1}A_{i+2}T^{-r+2i+3}T)$\\
\hspace*{1.65cm} $= \rho(A_{i+3}TA_{i+4}T\ldots A_{r-i-2}TA_{r-i-1}A_{i+2}T^{-r+2i+4})$\\
where $A_{r-i-1}A_{i+2}=A_{i+3}=A_{i+4}=\ldots =A_{r-i-2}$. \vspace{0.1cm}

Along the same way, we can prove also that if $r\geq 3$ is an odd integer and $x=T^{1-r}A_1TA_2T\ldots A_{r-1}TA_r \in E_1$, and under the same conditions on $A_i's$ mentioned in (a), (b) and (c), then there exists $w \in C_3$ such that $\rho(w^{-1}xw)=\rho(T^{-s}\underbrace{AT}_1AT\ldots \underbrace{AT}_s)$, where $A \in \{\alpha_1, \alpha_2, \alpha_1\alpha_2, \alpha_2\alpha_1, \alpha_1\alpha_2\alpha_1\}$ and $s\in \mathbb{N}$.
\end{proof}

\begin{example}
We consider three examples when $r=3$, $r=5$, and $r=7$.
\end{example}
\noindent (i) \underline{$r=3$}: Let $x=A_1TA_2TA_3T^{-2}$. \vspace{0.1cm} \\
If (a) holds true, that is $A_3A_1=A_2$, we take $w=A_1T$.\\
$\rho(w^{-1}xw)=\rho((A_1T)^{-1}A_1TA_2TA_3T^{-2}(A_1T))\\
\hspace*{1.55cm}= \rho(A_2TA_3T^{-2}A_1T) \\
\hspace*{1.55cm}= \rho(A_2)\rho(T)\rho(A_3)\rho(T^{-2})\rho(A_1)\rho(T) \\
\hspace*{1.55cm}= \rho(A_2)\rho(T)\rho(A_3)\rho(A_1)\rho(T^{-2})\rho(T) \\
\hspace*{1.55cm}= \rho(A_2TA_3A_1T^{-2}T)\\
\hspace*{1.55cm}= \rho(A_2TA_3A_1T^{-1})$, where $A_3A_1=A_2$. \vspace{0.1cm}\\
If (b) holds true, that is $A_3A_1=1$, take $w=A_1T$.\\
$ \rho(w^{-1}xw)= \rho((A_1T)^{-1}A_1TA_2TA_3T^{-2}(A_1T))\\
\hspace*{1.55cm}= \rho(A_2TA_3T^{-2}A_1T)\\
\hspace*{1.55cm}= \rho(A_2)\rho(T)\rho(A_3)\rho(T^{-2})\rho(A_1)\rho(T)\\
\hspace*{1.55cm}= \rho(A_2)\rho(T)\rho(A_3)\rho(A_1)\rho(T^{-2})\rho(T)\\
\hspace*{1.55cm}= \rho(A_2)\rho(T)\rho(A_3A_1)\rho(T^{-2})\rho(T)\\
\hspace*{1.55cm}= \rho(A_2)\rho(T)\rho(T^{-2})\rho(T)\\
\hspace*{1.55cm}=\rho(A_2TT^{-2}T)\\
\hspace*{1.55cm}= \rho(A_2)$.\vspace{0.1cm}\\
(ii) \underline{$r=5$}: Let $x=A_1TA_2TA_3TA_4TA_5T^{-4}$.\vspace{0.1cm}\\
If (a) holds true, that is $A_5A_1=A_2=A_3=A_4$, we take $w=A_1T$.\\
$\rho(w^{-1}xw)=\rho((A_1T)^{-1}A_1TA_2TA_3TA_4TA_5T^{-4}(A_1T))\\
\hspace*{1.55cm}= \rho(A_2TA_3TA_4TA_5T^{-4}A_1T)\\
\hspace*{1.55cm}= \rho(A_2)\rho(T)\rho(A_3)\rho(T)\rho(A_4)\rho(T)\rho(A_5)\rho(T^{-4})\rho(A_1)\rho(T)\\
\hspace*{1.55cm}= \rho(A_2)\rho(T)\rho(A_3)\rho(T)\rho(A_4)\rho(T)\rho(A_5)\rho(A_1)\rho(T^{-4})\rho(T)\\
\hspace*{1.55cm}= \rho(A_2TA_3TA_4TA_5A_1T^{-4}T)\\
\hspace*{1.55cm}= \rho(A_2TA_3TA_4TA_5A_1T^{-3})$, where $A_5A_1=A_2=A_3=A_4$. \vspace{0.1cm}\\ 
If (b) holds true, that is $A_5A_1=A_4A_2=1$, we take $w=A_1TA_2T$.\\
$\rho(w^{-1}xw)=\rho((A_1TA_2T)^{-1}A_1TA_2TA_3TA_4TA_5T^{-4}(A_1TA_2T))\\
\hspace*{1.55cm}=\rho(A_3TA_4TA_5T^{-4}A_1TA_2T)\\
\hspace*{1.55cm}=\rho(A_3)\rho(T)\rho(A_4)\rho(T)\rho(A_5)\rho(T^{-4})\rho(A_1)\rho(T)\rho(A_2)\rho(T)\\
\hspace*{1.55cm}=\rho(A_3)\rho(T)\rho(A_4)\rho(T)\rho(A_5)\rho(A_1)\rho(T^{-4})\rho(T)\rho(A_2)\rho(T)\\
\hspace*{1.55cm}=\rho(A_3)\rho(T)\rho(A_4)\rho(T)\rho(A_5A_1)\rho(T^{-4})\rho(T)\rho(A_2)\rho(T)\\
\hspace*{1.55cm}=\rho(A_3)\rho(T)\rho(A_4)\rho(T)\rho(T^{-4})\rho(T)\rho(A_2)\rho(T)\\
\hspace*{1.55cm}=\rho(A_3)\rho(T)\rho(A_4)\rho(T^{-2})\rho(A_2)\rho(T)\\
\hspace*{1.55cm}=\rho(A_3)\rho(T)\rho(A_4)\rho(A_2)\rho(T^{-2})\rho(T)\\
\hspace*{1.55cm}=\rho(A_3)\rho(T)\rho(A_4A_2)\rho(T^{-2})\rho(T)\\
\hspace*{1.55cm}=\rho(A_3)\rho(T)\rho(T^{-2})\rho(T)\\
\hspace*{1.55cm}= \rho(A_3TT^{-2}T)\\
\hspace*{1.55cm}=\rho(A_3)$.\\
If (c) holds true, we have just one case when $i=0$, that is $A_5A_1=1$ and $A_4A_2=A_3$. We take here $w_0=A_1TA_2T$.\\
$\rho(w_0^{-1}xw_0)=\rho((A_1TA_2T)^{-1}A_1TA_2TA_3TA_4TA_5T^{-4}(A_1TA_2T))\\
\hspace*{1.7cm}=\rho(A_3TA_4TA_5T^{-4}A_1TA_2T)\\
\hspace*{1.7cm}=\rho(A_3)\rho(T)\rho(A_4)\rho(T)\rho(A_5)\rho(T^{-4})\rho(A_1)\rho(T)\rho(A_2)\rho(T)\\
\hspace*{1.7cm}=\rho(A_3)\rho(T)\rho(A_4)\rho(T)\rho(A_5)\rho(A_1)\rho(T^{-4})\rho(T)\rho(A_2)\rho(T)\\
\hspace*{1.7cm}=\rho(A_3)\rho(T)\rho(A_4)\rho(T)\rho(A_5A_1)\rho(T^{-4})\rho(T)\rho(A_2)\rho(T)\\
\hspace*{1.7cm}=\rho(A_3)\rho(T)\rho(A_4)\rho(T)\rho(T^{-4})\rho(T)\rho(A_2)\rho(T)\\
\hspace*{1.7cm}=\rho(A_3)\rho(T)\rho(A_4)\rho(T^{-2})\rho(A_2)\rho(T)\\
\hspace*{1.7cm}=\rho(A_3)\rho(T)\rho(A_4)\rho(A_2)\rho(T^{-2})\rho(T)\\
\hspace*{1.08cm}=\rho(A_3TA_4A_2T^{-2}T),\\
\hspace*{1.08cm}=\rho(A_3TA_4A_2T^{-1})$, where $A_4A_2=A_3$.\vspace{0.1cm}\\
(iii) \underline{$r=7$}: Let $x=A_1TA_2TA_3TA_4TA_5TA_6TA_7T^{-6}$.\vspace{0.1cm}\\
If (a) holds true, that is $A_7A_1=A_2=A_3=A_4=A_5=A_6$, we take $w=A_1T$.\\
$\rho(w^{-1}xw)=\rho((A_1T)^{-1}A_1TA_2TA_3TA_4TA_5TA_6TA_7T^{-6}(A_1T))\\
\hspace*{1.55cm}= \rho(A_2TA_3TA_4TA_5TA_6TA_7T^{-6}A_1T)\\
\hspace*{1.55cm}= \rho(A_2)\rho(T)\rho(A_3)\rho(T)\rho(A_4)\rho(T)\rho(A_5)\rho(T)\rho(A_6)\rho(T)\rho(A_7)\rho(T^{-6})\rho(A_1)\rho(T)\\
\hspace*{1.55cm}= \rho(A_2)\rho(T)\rho(A_3)\rho(T)\rho(A_4)\rho(T)\rho(A_5)\rho(T)\rho(A_6)\rho(T)\rho(A_7)\rho(A_1)\rho(T^{-6})\rho(T)\\
\hspace*{1.55cm}= \rho(A_2TA_3TA_4TA_5TA_6TA_7A_1T^{-6}T)\\
\hspace*{1.55cm}= \rho(A_2TA_3TA_4TA_5TA_6TA_7A_1T^{-5}),$ where $A_7A_1=A_2=A_3=A_4=A_5=A_6$. \vspace{0.1cm}\\ 
If (b) holds true, that is $A_7A_1=A_6A_2=A_5A_3=1$, we take $w=A_1TA_2TA_3T$.\\
$\rho(w^{-1}xw)=\rho((A_1TA_2TA_3T)^{-1}A_1TA_2TA_3TA_4TA_5TA_6TA_7T^{-6}(A_1TA_2TA_3T))\\
\hspace*{1.55cm}= \rho(A_4TA_5TA_6TA_7T^{-6}A_1TA_2TA_3T)\\
\hspace*{1.55cm}= \rho(A_4)\rho(T)\rho(A_5)\rho(T)\rho(A_6)\rho(T)\rho(A_7)\rho(T^{-6})\rho(A_1)\rho(T)\rho(A_2)\rho(T)\rho(A_3)\rho(T)\\
\hspace*{1.55cm}= \rho(A_4)\rho(T)\rho(A_5)\rho(T)\rho(A_6)\rho(T)\rho(A_7)\rho(A_1)\rho(T^{-6})\rho(T)\rho(A_2)\rho(T)\rho(A_3)\rho(T)\\
\hspace*{1.55cm}= \rho(A_4)\rho(T)\rho(A_5)\rho(T)\rho(A_6)\rho(T)\rho(A_7A_1)\rho(T^{-6})\rho(T)\rho(A_2)\rho(T)\rho(A_3)\rho(T)\\
\hspace*{1.55cm}= \rho(A_4)\rho(T)\rho(A_5)\rho(T)\rho(A_6)\rho(T)\rho(T^{-6})\rho(T)\rho(A_2)\rho(T)\rho(A_3)\rho(T)\\
\hspace*{1.55cm}= \rho(A_4)\rho(T)\rho(A_5)\rho(T)\rho(A_6)\rho(T^{-4})\rho(A_2)\rho(T)\rho(A_3)\rho(T)\\
\hspace*{1.55cm}= \rho(A_4)\rho(T)\rho(A_5)\rho(T)\rho(A_6)\rho(A_2)\rho(T^{-4})\rho(T)\rho(A_3)\rho(T)\\
\hspace*{1.55cm}= \rho(A_4)\rho(T)\rho(A_5)\rho(T)\rho(A_6A_2)\rho(T^{-4})\rho(T)\rho(A_3)\rho(T)\\
\hspace*{1.55cm}= \rho(A_4)\rho(T)\rho(A_5)\rho(T)\rho(T^{-4})\rho(T)\rho(A_3)\rho(T)\\
\hspace*{1.55cm}= \rho(A_4)\rho(T)\rho(A_5)\rho(T^{-2})\rho(A_3)\rho(T)\\
\hspace*{1.55cm}= \rho(A_4)\rho(T)\rho(A_5)\rho(A_3)\rho(T^{-2})\rho(T)\\
\hspace*{1.55cm}= \rho(A_4)\rho(T)\rho(A_5A_3)\rho(T^{-2})\rho(T)\\
\hspace*{1.55cm}= \rho(A_4)\rho(T)\rho(T^{-2})\rho(T)\\
\hspace*{1.55cm}= \rho(A_4TT^{-2}T)\\
\hspace*{1.55cm}= \rho(A_4).\\ $
If (c) holds true, then we have 2 cases.\\
Case 1: $i=0$, then we have $A_7A_1=1$ and $A_6A_2=A_3=A_4=A_5$. We take here $w_0=A_1TA_2T$.\\
$\rho(w_0^{-1}xw_0)=\rho((A_1TA_2T)^{-1}A_1TA_2TA_3TA_4TA_5TA_6TA_7T^{-6}(A_1TA_2T))\\
\hspace*{1.70cm}=\rho(A_3TA_4TA_5TA_6TA_7T^{-6}A_1TA_2T)\\
\hspace*{1.70cm}=\rho(A_3)\rho(T)\rho(A_4)\rho(T)\rho(A_5)\rho(T)\rho(A_6)\rho(T)\rho(A_7)\rho(T^{-6})\rho(A_1)\rho(T)\rho(A_2)\rho(T)\\
\hspace*{1.70cm}=\rho(A_3)\rho(T)\rho(A_4)\rho(T)\rho(A_5)\rho(T)\rho(A_6)\rho(T)\rho(A_7)\rho(A_1)\rho(T^{-6})\rho(T)\rho(A_2)\rho(T)\\
\hspace*{1.70cm}=\rho(A_3)\rho(T)\rho(A_4)\rho(T)\rho(A_5)\rho(T)\rho(A_6)\rho(T)\rho(A_7A_1)\rho(T^{-6})\rho(T)\rho(A_2)\rho(T)\\
\hspace*{1.70cm}=\rho(A_3)\rho(T)\rho(A_4)\rho(T)\rho(A_5)\rho(T)\rho(A_6)\rho(T)\rho(T^{-6})\rho(T)\rho(A_2)\rho(T)\\
\hspace*{1.70cm}=\rho(A_3)\rho(T)\rho(A_4)\rho(T)\rho(A_5)\rho(T)\rho(A_6)\rho(T^{-4})\rho(A_2)\rho(T)\\
\hspace*{1.70cm}=\rho(A_3)\rho(T)\rho(A_4)\rho(T)\rho(A_5)\rho(T)\rho(A_6)\rho(A_2)\rho(T^{-4})\rho(T)\\
\hspace*{1.70cm}=\rho(A_3TA_4TA_5TA_6A_2T^{-4}T)\\
\hspace*{1.70cm}=\rho(A_3TA_4TA_5TA_6A_2T^{-3}),$ where $A_6A_2=A_3=A_4=A_5$.\vspace{0.1cm}\\
Case 2: $i=1$, then we have $A_7A_1=A_6A_2=1$ and $A_5A_3=A_4$. We take here $w_1=A_1TA_2TA_3T$.\\
$\rho(w_1^{-1}xw_1)=\rho((A_1TA_2TA_3T)^{-1}A_1TA_2TA_3TA_4TA_5TA_6TA_7T^{-6}(A_1TA_2TA_3T))\\
\hspace*{1.70cm}=\rho(A_4TA_5TA_6TA_7T^{-6}A_1TA_2TA_3T)\\
\hspace*{1.70cm}=\rho(A_4)\rho(T)\rho(A_5)\rho(T)\rho(A_6)\rho(T)\rho(A_7)\rho(T^{-6})\rho(A_1)\rho(T)\rho(A_2)\rho(T)\rho(A_3)\rho(T)\\
\hspace*{1.70cm}=\rho(A_4)\rho(T)\rho(A_5)\rho(T)\rho(A_6)\rho(T)\rho(A_7)\rho(A_1)\rho(T^{-6})\rho(T)\rho(A_2)\rho(T)\rho(A_3)\rho(T)\\
\hspace*{1.70cm}=\rho(A_4)\rho(T)\rho(A_5)\rho(T)\rho(A_6)\rho(T)\rho(A_7A_1)\rho(T^{-6})\rho(T)\rho(A_2)\rho(T)\rho(A_3)\rho(T)\\
\hspace*{1.70cm}=\rho(A_4)\rho(T)\rho(A_5)\rho(T)\rho(A_6)\rho(T)\rho(T^{-6})\rho(T)\rho(A_2)\rho(T)\rho(A_3)\rho(T)\\
\hspace*{1.70cm}=\rho(A_4)\rho(T)\rho(A_5)\rho(T)\rho(A_6)\rho(T^{-4})\rho(A_2)\rho(T)\rho(A_3)\rho(T)\\
\hspace*{1.70cm}=\rho(A_4)\rho(T)\rho(A_5)\rho(T)\rho(A_6)\rho(A_2)\rho(T^{-4})\rho(T)\rho(A_3)\rho(T)\\
\hspace*{1.70cm}=\rho(A_4)\rho(T)\rho(A_5)\rho(T)\rho(A_6A_2)\rho(T^{-4})\rho(T)\rho(A_3)\rho(T)\\
\hspace*{1.70cm}=\rho(A_4)\rho(T)\rho(A_5)\rho(T)\rho(T^{-4})\rho(T)\rho(A_3)\rho(T)\\
\hspace*{1.70cm}=\rho(A_4)\rho(T)\rho(A_5)\rho(T^{-2})\rho(A_3)\rho(T)\\
\hspace*{1.70cm}=\rho(A_4)\rho(T)\rho(A_5)\rho(A_3)\rho(T^{-2})\rho(T)\\
\hspace*{1.70cm}=\rho(A_4TA_5A_3T^{-2}T)\\
\hspace*{1.70cm}=\rho(A_4TA_5A_3T^{-1}), $ where $A_5A_3=A_4$.

\begin{theorem}
Suppose $q^{6k}\neq 1$ for all $k \in \mathbb{Z}$, $q \notin P_m\cup R_m$ for any even integer $m$ and $q\notin S_n$ for any  integer $n$. Let $r\geq 3$ be an odd integer and let $x\in E_1$. If $x$ satisfies the hypothesis of Proposition 9, then $x \notin \ker \rho $.
\end{theorem}
\begin{proof}
If $x$ satisfies the hypothesis of Proposition 9, then, without loss of generality, there exists  $w \in C_3$ such that $\rho(w^{-1}xw)=\rho(ATAT\ldots \underbrace{AT}_sT^{-s})$, where $s\in \mathbb{N}$. If $x \in \ker \rho$, then $\rho(x)=I_3$, and so $\rho(ATAT\ldots \underbrace{AT}_sT^{-s})=\rho(w^{-1}xw)=\rho(w^{-1})\rho(x)\rho(w)=\rho(w^{-1})\rho(w)=\rho(w^{-1}w)=I_3$, which contradicts Theorem 8. So $x \notin \ker \rho$.
\end{proof}

\vskip .2in

{\bf Conflict of Interest.} On behalf of all authors, the corresponding author states that there is no conflict of interest.

\end{document}